# Surface family with a common natural asymptotic lift of a timelike curve in Minkowski 3-space

[1] Ergin Bayram, [2] Evren Ergün, [1] Emin Kasap

[1] Ondokuz Mayıs University, Faculty of Arts and Sciences, Department of Mathematics, Samsun, Turkey.

[2] Ondokuz Mayıs University, Çarşamba Chamber of Commerce Vocational School, Çarşamba, Samsun, Turkey.

## Abstract

In the present paper, we find a surface family possessing the natural lift of a given timelike curve as a asymptotic in Minkowski 3-space. We express necessary and sufficient conditions for the given curve such that its natural lift is a asymptotic on any member of the surface family. Finally, we illustrate the method with some examples.

## Introduction and Preliminaries

The problem of finding surfaces with a given common curve as a special curve was firstly handled by Wang et.al. [1]. They constructed surfaces with a common geodesic. Li et.al. [2] obtained necessary and sufficient condition for a given curve to be a line of curvature on a surface pencil. Bayram et.al. [3] studied surface pencil with a common asymptotic curve. In 2014, Ergün et.al. presented a constraint for surfaces with a common line of curvature in Minkowski 3-space. In [4-9] authors found necessary and sufficient conditions for a given curve to be a common special curve on a surface family in Euclidean and Minkowski spaces.

Inspired with the above studies, we find necessary and sufficient condition for a surface family possessing the natural lift of a given timelike curve as a common asymptotic curve. First, we start with fundamentals required for the paper.

Minkowski 3-space $\mathbb{R}_1^3$ is the vector space $\mathbb{R}^3$ equipped with the Lorentzian inner product g given by

$$g(X,X) = -x_1^2 + x_2^2 + x_3^2$$

where $X = (x_1, x_2, x_3) \in \mathbb{R}^3$. A vector $X = (x_1, x_2, x_3) \in \mathbb{R}^3$ is said to be timelike if $g(X,X) < 0$, spacelike if $g(X,X) > 0$ or $X = 0$ and lightlike (or null) if $g(X,X) = 0$ and $x \neq 0$ [10]. Similarly, an arbitrary curve $\alpha = \alpha(s)$ in $\mathbb{R}_1^3$ can locally be timelike, spacelike or null (lightlike), if all of its velocity vectors $\alpha'(s)$ are respectively timelike, spacelike or null (lightlike), for every $s \in I \subset \mathbb{R}$. A lightlike vector $X$ is said to be positive (resp. negative) if and only if $x_1 > 0$ (resp. $x_1 < 0$) and a timelike vector $X$ is said to be positive (resp. negative) if and only if $x_1 > 0$ (resp. $x_1 < 0$). The norm of a vector $X$ is defined by $\|X\|_{IL} = \sqrt{|g(X,X)|}$ [10].

The vectors $X = (x_1, x_2, x_3)$, $Y = (y_1, y_2, y_3) \in \mathbb{R}^3_1$ are orthogonal if and only if $g(X,X) = 0$ [11].

Now let $X$ and $Y$ be two vectors in $\mathbb{R}^3_1$, then the Lorentzian cross product is given by [12]

$$X \times Y = \begin{vmatrix} \vec{e}_1 & -\vec{e}_2 & -\vec{e}_3 \\ x_1 & x_2 & x_3 \\ y_1 & y_2 & y_3 \end{vmatrix}$$
$$= (x_2 y_3 - x_3 y_2,\ x_1 y_3 - x_3 y_1,\ x_2 y_1 - x_1 y_2).$$

We denote by $\{T(s), N(s), B(s)\}$ the moving Frenet frame along the curve $\alpha$. Then $T$, $N$ and $B$ are the tangent, the principal normal and the binormal vector of the curve $\alpha$, respectively.

Let $\alpha$ be a unit speed timelike curve with curvature $\kappa$ and torsion $\tau$. So, $T$ is a timelike vector field, $N$ and $B$ are spacelike vector fields. For these vectors, we can write

$$T \times N = -B,\quad N \times B = T,\quad B \times T = -N,$$

where $\times$ is the Lorentzian cross product in $\mathbb{R}^3_1$ [13]. The binormal vector field $B(s)$ is the unique spacelike unit vector field perpendicular to the timelike plane $\{T(s), N(s)\}$ at every point $\alpha(s)$ of $\alpha$, such that $\{T, N, B\}$ has the same orientation as $\mathbb{R}^3_1$. Then, Frenet formulas are given by [13]

$$T' = \kappa N,\ N' = \kappa T + \tau B,\ B' = -\tau N.$$

Let $\alpha$ be a unit speed spacelike curve with spacelike binormal. Now, $T$ and $B$ are spacelike vector fields and $N$ is a timelike vector field. In this situation,

$$T \times N = -B,\quad N \times B = -T,\quad B \times T = N.$$

The binormal vector field $B(s)$ is the unique spacelike unit vector field perpendicular to the timelike plane $\{T(s), N(s)\}$ at every point $\alpha(s)$ of $\alpha$, such that $\{T, N, B\}$ has the same orientation as $\mathbb{R}^3_1$. Then, Frenet formulas are given by [13]

$$T' = \kappa N,\ N' = \kappa T + \tau B,\ B' = \tau N.$$

Let $\alpha$ be a unit speed spacelike curve with timelike binormal. Now, $T$ and $N$ are spacelike vector fields and $B$ is a timelike vector field. In this situation,

$$T \times N = B, \quad N \times B = -T, \quad B \times T = -N,$$

The binormal vector field $B(s)$ is the unique timelike unit vector field perpendicular to the spacelike plane $\{T(s), N(s)\}$ at every point $\alpha(s)$ of $\alpha$, such that $\{T, N, B\}$ has the same orientation as $\mathbb{R}_1^3$. Then, Frenet formulas are given by [13]

$$T' = \kappa N, \ N' = -\kappa T + \tau B, \ B' = \tau N.$$

**Lemma** *Let $X$ and $Y$ be nonzero Lorentz orthogonal vectors in $\mathbb{R}_1^3$. If $X$ is timelike, then $Y$ is spacelike [11].*

Lemma *Let $X$ and $Y$ be positive (negative) timelike vectors in $\mathbb{R}_1^3$. Then*

$$g(X, Y) \leq \|X\| \|Y\|$$

*whit equality if and only if $X$ and $Y$ are linearly dependent [11].*

**Lemma**

*i) Let $X$ and $Y$ be positive (negative) timelike vectors in $\mathbb{R}_1^3$. By Lemma 2, there is a unique nonnegative real number $\varphi(X, Y)$ such that*

$$g(X, Y) = \|X\| \|Y\| \cosh \varphi(X, Y).$$

*The Lorentzian timelike angle between $X$ and $Y$ is defined to be $\varphi(X, Y)$ [11].*

*ii) Let $X$ and $Y$ be spacelike vectors in $\mathbb{R}_1^3$ that span a spacelike vector subspace. Then we have*

$|g(X, Y)| \leq \|X\| \|Y\|.$

*Hence, there is a unique real number $\varphi(X, Y)$ between $0$ and $\pi$ such that*

$g(X, Y) = \|X\| \|Y\| \cos \varphi(X, Y).$

*$\varphi(X, Y)$ is defined to be the Lorentzian spacelike angle between $X$ and $Y$ [11].*

*iii) Let $X$ and $Y$ be spacelike vectors in $\mathbb{R}_1^3$ that span a timelike vector subspace. Then, we have*

$g(X, Y) > \|X\| \|Y\|.$

*Hence, there is a unique positive real number $\varphi(X, Y)$ between $0$ and $\pi$ such that*

$|g(X, Y)| = \|X\| \|Y\| \cosh \varphi(X, Y).$

*$\varphi(X, Y)$ is defined to be the Lorentzian timelike angle between $X$ and $Y$ [11].*

***iv)*** *Let $X$ be a spacelike vector and $Y$ be a positive timelike vector in $\mathbb{R}_1^3$. Then there is a unique nonnegative real number $\varphi(X,Y)$ such that*

$$|g(X,Y)| = \|X\|\|Y\|\sinh\varphi(X,Y).$$

*$\varphi(X,Y)$ is defined to be the Lorentzian timelike angle between $X$ and $Y$ [11].*

Let $P$ be a surface in $\mathbb{R}_1^3$ and let $\alpha : I \to P$ be a parametrized curve. $\alpha$ is called an integral curve of $X$ if

$$\frac{d}{ds}(\alpha(s)) = X(\alpha(s)), \; \left(\text{for all } t \in I\right),$$

where $X$ is a smooth tangent vector field on $P$ [10]. We have

$$TP = \bigcup_{p \in P} T_p P = \chi(P),$$

where $T_p P$ is the tangent space of $P$ at $p$ and $\chi(P)$ is the space of tangent vector fields on $P$.

For any parametrized curve $\alpha : I \to P$, $\bar{\alpha} : I \to TP$ is given by

$$\bar{\alpha}(s) = (\alpha(s), \alpha'(s)) = \alpha'(s)|_{\alpha(s)}$$

is called the natural lift of $\alpha$ on $TP$ [14].

Let $\alpha(s)$, $L_1 \leq s \leq L_2$, be an arc length timelike curve. Then, the natural lift $\bar{\alpha}$ of $\alpha$ is a spacelike curve with timelike or spacelike binormal. We have following relations between the Frenet frame $\{T(s), N(s), B(s)\}$ of $\alpha$ and the Frenet frame $\{\bar{T}(s), \bar{N}(s), \bar{B}(s)\}$ of $\bar{\alpha}$.

**a)** Let the natural lift $\bar{\alpha}$ of $\alpha$ is a spacelike curve with timelike binormal.

**i)** If the Darboux vector $W$ of the curve $\alpha$ is a spacelike vector, then we have

$$\begin{pmatrix} \bar{T}(s) \\ \bar{N}(s) \\ \bar{B}(s) \end{pmatrix} = \begin{pmatrix} 0 & 1 & 0 \\ \cosh\theta & 0 & \sinh\theta \\ \sinh\theta & 0 & \cosh\theta \end{pmatrix} \begin{pmatrix} T(s) \\ N(s) \\ B(s) \end{pmatrix}. \tag{1}$$

**ii)** If $W$ is a timelike vector, then we have

$$\begin{pmatrix} \bar{T}(s) \\ \bar{N}(s) \\ \bar{B}(s) \end{pmatrix} = \begin{pmatrix} 0 & 1 & 0 \\ \sinh\theta & 0 & \cosh\theta \\ \cosh\theta & 0 & \sinh\theta \end{pmatrix} \begin{pmatrix} T(s) \\ N(s) \\ B(s) \end{pmatrix}. \tag{2}$$

**b)** Let the natural lift $\bar{\alpha}$ of $\alpha$ is a spacelike curve with spacelike binormal.

**i)** If $W$ is a spacelike vector, then we have

$$\begin{pmatrix} \bar{T}(s) \\ \bar{N}(s) \\ \bar{B}(s) \end{pmatrix} = \begin{pmatrix} 0 & 1 & 0 \\ \cosh\theta & 0 & \sinh\theta \\ -\sinh\theta & 0 & -\cosh\theta \end{pmatrix} \begin{pmatrix} T(s) \\ N(s) \\ B(s) \end{pmatrix}. \tag{3}$$

**ii)** If $W$ is a timelike vector, then we have

$$\begin{pmatrix} \bar{T}(s) \\ \bar{N}(s) \\ \bar{B}(s) \end{pmatrix} = \begin{pmatrix} 0 & 1 & 0 \\ \sinh\theta & 0 & \cosh\theta \\ -\cosh\theta & 0 & -\sinh\theta \end{pmatrix} \begin{pmatrix} T(s) \\ N(s) \\ B(s) \end{pmatrix}. \tag{4}$$

**Surface family with a common natural asymptotic lift of a timelike curve in Minkowski 3-space**

Suppose we are given a 3-dimensional timelike curve $\alpha(s)$, $L_1 \leq s \leq L_2$, in which $s$ is the arc length and $\|\alpha''(s)\| \neq 0$, $L_1 \leq s \leq L_2$. Let $\bar{\alpha}(s)$, $L_1 \leq s \leq L_2$, be the natural lift of the given curve $\alpha(s)$. Now, $\bar{\alpha}$ is a spacelike curve with timelike or spacelike binormal. Surface family that interpolates $\bar{\alpha}(s)$ as a common curve is given in the parametric form as

$$P(s,t) = \bar{\alpha}(s) + u(s,t)\bar{T}(s) + v(s,t)\bar{N}(s) + w(s,t)\bar{B}(s), \tag{5}$$

where $u(s,t)$, $v(s,t)$ and $w(s,t)$ are $C^1$ functions, called *marching-scale functions,* and $\{\bar{T}(s), \bar{N}(s), \bar{B}(s)\}$ is the Frenet frame of the curve $\bar{\alpha}$.

*Remark Observe that choosing different marching-scale functions yields different surfaces possessing $\bar{\alpha}(s)$ as a common curve.*

Our goal is to find the necessary and sufficient conditions for which the curve $\bar{\alpha}(s)$ is isoparametric and asymptotic curve on the surface $P(s,t)$. Firstly, as $\bar{\alpha}(s)$ is an isoparametric curve on the surface $P(s,t)$, there exists a parameter $t_0 \in [T_1, T_2]$ such that

$$u(s, t_0) = v(s, t_0) = w(s, t_0) \equiv 0, \ L_1 \leq s \leq L_2, \ T_1 \leq t_0 \leq T_2. \tag{6}$$

Secondly the curve $\bar{\alpha}$ is an asymptotic curve on the surface $P(s,t)$ if and only if along the curve the normal vector field $n(s,t_0)$ of the surface is parallel to the binormal vector field $\bar{B}$ of the curve $\bar{\alpha}$. The normal vector of $P(s,t)$ can be written as

$$n(s,t) = \frac{\partial P(s,t)}{\partial s} \times \frac{\partial P(s,t)}{\partial t}.$$

Along the curve $\bar{\alpha}$, one can obtain the normal vector $n(s,t_0)$ using Eqns. (ref. 5 – ref. 6) with an appropriate equation in Eqns. (1-4). It has one of the following forms:

**i)** if $\bar{\alpha}$ is a spacelike curve with timelike binormal and Darboux vector $W$ is spacelike or timelike, then

$$n(s,t_0) = \kappa\left[\frac{\partial w}{\partial t}(s,t_0)\bar{N}(s) + \frac{\partial v}{\partial t}(s,t_0)\bar{B}(s)\right], \tag{7}$$

**ii)** if $\bar{\alpha}$ is a spacelike curve with spacelike binormal and Darboux vector $W$ is spacelike, then

$$n(s,t_0) = -\kappa\left[\frac{\partial w}{\partial t}(s,t_0)\bar{N}(s) + \frac{\partial v}{\partial t}(s,t_0)\bar{B}(s)\right], \tag{8}$$

**iii)** if $\bar{\alpha}$ is a spacelike curve with spacelike binormal and Darboux vector $W$ timelike, then

$$n(s,t_0) = \kappa\left[\frac{\partial w}{\partial t}(s,t_0)\bar{N}(s) - \frac{\partial v}{\partial t}(s,t_0)\bar{B}(s)\right], \tag{9}$$

where $\kappa$ is the curvature of the curve $\alpha$.

Since $\kappa(s) \neq 0$, $L_1 \leq s \leq L_2$, the curve $\bar{\alpha}$ is an asymptotic curve on the surface $P(s,t)$ if and only if

$$\frac{\partial w}{\partial t}(s,t_0) \equiv 0, \quad \frac{\partial v}{\partial t}(s,t_0) \neq 0.$$

So, we can present :

**Theorem** Let $\alpha(s)$, $L_1 \leq s \leq L_2$, be a unit speed timelike curve with nonvanishing curvature and $\bar{\alpha}(s)$, $L_1 \leq s \leq L_2$, be its natural lift. $\bar{\alpha}$ is an asymptotic curve on the surface (5) if and only if

$$\begin{cases} u(s,t_0) = v(s,t_0) = w(s,t_0) = \frac{\partial w}{\partial t}(s,t_0) \equiv 0, \\ \frac{\partial v}{\partial t}(s,t_0) \neq 0, \end{cases} \tag{10}$$

where $L_1 \leq s \leq L_2$, $T_1 \leq t$, $t_0 \leq T_2$ ($t_0$ fixed).

**Corollary** Let $\alpha(s)$, $L_1 \leq s \leq L_2$, be a unit speed timelike curve with nonvanishing curvature and $\bar{\alpha}(s)$, $L_1 \leq s \leq L_2$, be its natural lift. If

$$u(s,t) = v(s,t) = t - t_0, \quad w(s,t) \equiv 0$$

or $\qquad (11)$

$$u(s,t) = w(s,t) \equiv 0, \quad v(s,t) = t - t_0,$$

where $L_1 \leq s \leq L_2$, $T_1 \leq t$, $t_0 \leq T_2$ ($t_0$ fixed) then (5) is a ruled surface possessing $\bar{\alpha}$

*as an asymptotic curve.*

Proof By taking marching scale functions as $u(s,t) = v(s,t) = t - t_0$, $w(s,t) \equiv 0$ or $u(s,t) = w(s,t) \equiv 0$, $v(s,t) = t - t_0$, the surface (5) takes the form

$$P(s,t) = \bar{\alpha}(s) + (t - t_0)[\bar{T}(s) + \bar{N}(s)]$$

or

$$P(s,t) = \bar{\alpha}(s) + (t - t_0)\bar{N}(s),$$

which is a ruled surface satisfying Eqn. (10).

**Examples**
**Example 1**

Let $\alpha(s) = (\sinh s, 0, \cosh s)$ be a timelike curve. It is easy to show that

$$T(s) = (\cosh s, 0, \sinh s),$$
$$N(s) = (\sinh s, 0, \cosh s),$$
$$B(s) = (0, -1, 0).$$

The natural lift of the curve $\alpha$ is $\bar{\alpha}(s) = (\cosh s, 0, \sinh s)$ and its Frenet vectors

$$\bar{T}(s) = (\sinh s, 0, \cosh s),$$
$$\bar{N}(s) = (\cosh s, 0, \sinh s),$$
$$\bar{B}(s) = (0, 1, 0).$$

Choosing marching scale functions as $u(s,t) = t$, $v(s,t) = \sinh t$, $w(s,t) \equiv 0$ Eqn. (10) is satisfied and we obtain the surface

$$P_1(s,t) = (\cosh s + t \sinh s + (\sinh t)\cosh s, \sinh t, t \cosh s + \sinh s + (\sinh t)\sinh s).$$

$-1 \leq s \leq 1$, $-1 \leq t \leq 0$, possessing $\bar{\alpha}$ as a common natural asymptotic lift (Fig. 1).
For the same curve, if we choose $u(s,t) \equiv 0$, $v(s,t) = (\sinh s)\sinh t$, $w(s,t) = t - \sinh t$ we get the surface

$$P_2(s,t) = ((\cosh s)(1 + (\sinh s)\sinh t), (t - \sinh t), (\sinh s)(1 + (\sinh s)\sinh t)),$$

$0 < s \leq 1$, $-1 \leq t \leq 1$, satisfying Eqn. (10) and accepting $\bar{\alpha}$ as a common natural asymptotic lift (Fig. 2).

**Example 2**

Given the arclength timelike curve $\alpha(s) = \left(\frac{5}{3}s, \frac{4}{9}\cos 3s, \frac{4}{9}\sin 3s\right)$ its Frenet apparatus are

$$T(s) = \left(\frac{5}{3}, -\frac{4}{3}\sin 3s, \frac{4}{3}\cos 3s\right),$$

$$N(s) = (0, -\cos 3s, -\sin 3s),$$

$$B(s) = \left(-\frac{4}{3}, \frac{5}{3}\sin 3s, -\frac{5}{3}\cos 3s\right).$$

The natural lift of the curve $\alpha$ is $\bar{\alpha}(s) = \left(\frac{5}{3}, -\frac{4}{3}\sin 3s, \frac{4}{3}\cos 3s\right)$ and its Frenet vectors

$$\bar{T}(s) = (0, -\cos 3s, -\sin 3s),$$

$$\bar{N}(s) = (0, \sin 3s, -\cos 3s),$$

$$\bar{B}(s) = (-1, 0, 0).$$

If we let marching scale functions as $u(s,t) = w(s,t) \equiv 0$, $v(s,t) = t$, we get the ruled surface

$$P_3(s,t) = \left(\frac{5}{3}, \left(t - \frac{4}{3}\right)\sin 3s, \left(\frac{4}{3} - t\right)\cos 3s\right),$$

$-1.1 \leq s \leq 1, -1 \leq t \leq 1$, satisfying Eqn. (11) and passing through $\bar{\alpha}$ as a common natural asymptotic lift (Fig. 3).

For the same curve, if we choose $u(s,t) \equiv 0$, $v(s,t) = t\ln s$, $w(s,t) = t^2 e^s$ we obtain the surface

$$P_4(s,t) = \left(\frac{5}{3} - t^2 e^s, \left(t\ln s - \frac{4}{3}\right)\sin(3s), \left(\frac{4}{3} - t\ln s\right)\cos(3s)\right),$$

$1 < s \leq 2, 0 \leq t \leq 1$, satisfying Eqn. (10) and possessing $\bar{\alpha}$ as a common natural asymptotic lift (Fig. 4).

**Acknowledgements**

The first author would like to thank TUBITAK (The Scientific and Technological Research Council of Turkey) for their financial supports during his doctorate studies.